# A HSS Matrix-Inspired Butterfly-Based Direct Solver for Analyzing Scattering from Two-dimensional Objects

Yang Liu[†], Han Guo[†], and Eric Michielssen[†]

**Abstract**

A butterfly-based fast direct integral equation solver for analyzing high-frequency scattering from two-dimensional objects is presented. The solver leverages a randomized butterfly scheme to compress blocks corresponding to near- and far-field interactions in the discretized forward and inverse electric field integral operators. The observed memory requirements and computational cost of the proposed solver scale as $O(N \log^2 N)$ and $O(N^{1.5} \log N)$, respectively. The solver is applied to the analysis of scattering from electrically large objects spanning over ten thousand of wavelengths and modeled in terms of five million unknowns.

## 1 Introduction

Fast direct integral equation (IE) solvers constitute attractive alternatives to fast multipole method-based iterative solvers due to their ability to solve electromagnetic scattering problems that are inherently ill-conditioned and/or involve many right-hand sides (RHS). Present direct solvers leverage the low-rank (LR) nature of off-diagonal blocks of the discretized IE operator or its inverse [1-4]. This property leads to direct solvers with quasi-linear CPU and memory requirements for electrically small and structured scatterers [2, 4]. However, for electrically large and arbitrarily-shaped scatterers, the CPU and memory requirements of these LR solvers deteriorate to $O(N^\alpha \log^\beta N)$ ($\alpha \geq 2$, $\beta \geq 1$) and $O(N^\alpha \log N)$ ($\alpha \geq 1.5$) as the blocks of the inverse are no longer LR compressible.

Recently, we developed a new class of direct solvers that leverage butterfly schemes, a.k.a. multilevel matrix decomposition algorithms (MLMDA), to compress blocks in the discretized IE operator and its hierarchical LU factorization [5]. Butterfly schemes [6, 7] leverage the LR nature of judiciously selected sub-blocks of off-diagonal blocks (that themselves are LR-incompressible) of a discretized forward operator and its LU factors. The resulting butterfly-based direct solvers attain $O(N \log^2 N)$ and $O(N^{1.5} \log N)$ memory and CPU complexities, irrespective of the nature and electrical size of the scatterer. These solvers have been successfully applied to the analyses of scattering from perfectly electrically conducting (PEC) and homogenous dielectric, two and three dimensional (2D and 3D) objects involving many million unknowns. That said, in the context of analyzing low frequency scattering phenomena, hierarchical LU factorizations are computationally expensive relative to recently developed techniques that leverage hierarchically semi-separable (HSS) matrices and related constructs [3, 4, 8] due to their simplicity and ability to compress blocks representing near-field interactions. Unfortunately, HSS matrix techniques developed to date only attain quasi-linear complexities for electrically small 2D scattering problems [3, 8].

[†]Department of Electrical Engineering and Computer Science, University of Michigan, Ann Arbor, MI (liuyangz@umich.edu, hanguo@umich.edu, emichiel@umich.edu).

This work develops a butterfly-based direct solver inspired by HSS matrix techniques for analyzing scattering from electrically large 2D objects. Specifically, the proposed solver factorizes the impedance matrix as a product of sparse factors; each factor consists of the identity matrix and butterfly-compressed off-diagonal blocks representing a (partial) scattering matrix involving adjacent subscatterers. The factorization and inversion process hinges on a fast randomized scheme capable of constructing an arbitrary-level butterfly factorization of a partial scattering matrix with overwhelmingly high probability. The proposed butterfly-based direct solver is applied to the analysis of scattering from electrically large objects spanning over ten thousand wavelengths.

## 2 Butterfly Compression of the Impedance Matrix

Let $S$ denote a PEC cylindrical shell residing in free space. A $TM_z$ field $E_z^{inc}(\boldsymbol{\rho})$ impinges on $S$ and induces a current $J_z(\boldsymbol{\rho})$. Requiring the total electric field on $S$ to vanish yields the electric field integral equation (EFIE)

$$E_z^{inc}(\boldsymbol{\rho}) = \frac{k\eta_0}{4} \int_S J_z(\boldsymbol{\rho}') H_0^{(2)}(k|\boldsymbol{\rho}-\boldsymbol{\rho}'|) ds' \quad \forall \boldsymbol{\rho} \in S. \tag{1}$$

Here, $k = 2\pi/\lambda_0$ is the wavenumber, $\lambda_0$ denotes the free-space wavelength, $\eta_0$ is the intrinsic impedance of free space, and $H_0^{(2)}(\cdot)$ is the zeroth-order Hankel function of the second kind. To numerically solve (1), $J_z(\boldsymbol{\rho})$ is expanded into $N$ basis functions as $J_z(\boldsymbol{\rho}) = \sum_{i=1}^{N} I_i b_i(\boldsymbol{\rho})$. Testing (1) with functions $t_i(\boldsymbol{\rho})$, $i=1,...,N$ yields the matrix equation

$$\boldsymbol{Z} \cdot \boldsymbol{I} = \boldsymbol{V}. \tag{2}$$

In (2) vector $\boldsymbol{I}$ collects the current expansion coefficients $I_j$, and the elements of $\boldsymbol{Z}$ and $\boldsymbol{V}$ are

$$\boldsymbol{Z}_{ij} = \frac{k\eta_0}{4} \int_S t_i(\boldsymbol{\rho}) \int_S b_j(\boldsymbol{\rho}') H_0^{(2)}(k|\boldsymbol{\rho}-\boldsymbol{\rho}'|) ds' ds \tag{3}$$

$$V_i = \int_S t_i(\boldsymbol{\rho}) E_z^{inc}(\boldsymbol{\rho}) ds. \tag{4}$$

The proposed direct solver starts from a butterfly-compressed approximation to $\boldsymbol{Z}$ constructed as follows.

First, $S$ is decomposed into two equal-sized level-1 subscatterers, each containing approximately $N/2$ basis functions. This step is repeated $L-1$ times, resulting in a binary tree with $L$ levels. At level $0 \leq l \leq L$, there are $2^l$ (sub)scatterers, each containing approximately $N/2^l$ basis functions.

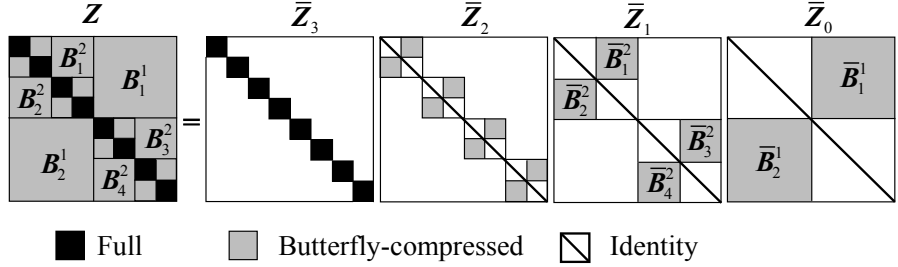

Figure 1: Butterfly-compression and factorization of the impedance matrix $Z$ with $L = 3$.

In what follows, submatrices of $Z$ that model self-interactions of level-$L$ subscatterers are denoted $Z_k^L$, $k = 1,...,2^L$ and directly computed using (3). Submatrices that model interactions between the two (adjacent) children of the $k^{th}$ level-$l$ subscatterer are denoted $B_{2k-1}^{l+1}$ and $B_{2k}^{l+1}$, $k = 1,...,2^l$, $l = 0,...,L-1$, and these submatrices are compressed by the butterfly scheme (Fig. 1). Specifically, the $m \times n$ submatrix $B_k^l$ with $m \approx n \approx N/2^l$ is recursively partitioned using $V = L - l$ levels by repeating the above-described process: at level $v = 0,...,V$, there are $2^v$ observation subscatterers containing approximately $m/2^v$ testing functions and $2^{V-v}$ source subscatterers containing approximately $n/2^{V-v}$ basis functions. Define the butterfly rank $r$ as the maximum (numerical) interaction rank for all $2^V(V+1)$ subscatterer pairs. Upon constructing the LR factorizations for all these subscatterer pairs, a $V$-level butterfly representation for submatrix $B_k^l$ is

$$B_k^l = R^{V+1} R^V \cdots R^1 R^0 \quad (5)$$

where $R^{V+1} = \text{diag}(R^{V+1,1},...,R^{V+1,2^V})$ and $R^0 = \text{diag}(R^{0,1},...,R^{0,2^V})$ are projection matrices. Their diagonal blocks $R^{V+1,i}$ and $R^{0,i}$ have approximate dimensions $(m/2^V) \times r$ and $r \times (n/2^V)$, respectively. The kernel matrices $R^v$, $v = 1,...,V$, consist of blocks of approximate dimensions $r \times 2r$ and are block diagonal following a row permutation, i.e., $D^v R^v = \text{diag}(R^{v,1},...,R^{v,2^{V-1}})$ where $D^v$ is the permutation matrix that renders $R^v$ diagonal and the diagonal blocks $R^{v,i}$ have approximate dimensions $2r \times 2r$ (Fig. 2).

It can be shown that the maximum butterfly rank $r$ for all submatrices $B_k^l$, $k = 1,...,2^l$, $l = 1,...,L$ stay in essence as constant. Note that for structures with sharp corners, it is required that each corner is fully contained in some level-$L$ subscatterer. As was shown in [6], the CPU and memory requirements for butterfly compressing the entire $Z$ matrix scale as $O(N \log^2 N)$.

## 3 Butterfly-based Inversion of the Impedance Matrix

Next, the impedance matrix $Z$ is factorized and inverted via the extension of a HSS-type LR solvers to butterflies [8]. Specifically, $Z$ is factorized as

$$Z = \bar{Z}_L \bar{Z}_{L-1} \cdots \bar{Z}_0 \quad (1)$$

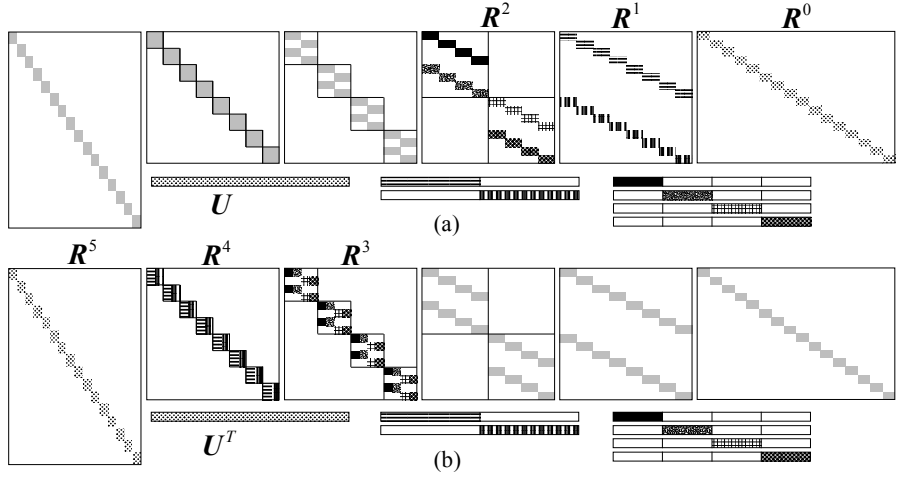

Figure 2: Random matrices $U$ and their associated blocks in (a) $R^{v_m},...,R^0$ and (b) $R^V,...,R^{v_m+1}$ for construction of a $L = 4$-level butterfly.

where each factor $\bar{Z}_l$, $l = 0,...,L$ is block diagonal as $\bar{Z}_l = \text{diag}(\bar{Z}_1^l,...,\bar{Z}_{2^l}^l)$ (Fig. 1). The diagonal blocks are $\bar{Z}_k^L = Z_k^L$, $k = 1,...,2^L$ and

$$\bar{Z}_k^{l-1} = \begin{bmatrix} I & \bar{B}_{2k-1}^l \\ \bar{B}_{2k}^l & I \end{bmatrix}, k = 1,...,2^l, l \leq L. \quad (2)$$

Here, $I$ denotes the identity matrix and the $\bar{B}_k^l$, $k = 1,...,2^l$, $l = 1,...,L$ are

$$\bar{B}_k^l = [\bar{Z}_k^l]^{-1} \begin{bmatrix} \bar{Z}_{2k-1}^{l+1} & \\ & \bar{Z}_{2k}^{l+1} \end{bmatrix}^{-1} ... \begin{bmatrix} \bar{Z}_{2^{L-l}(k-1)+1}^L & & \\ & \ddots & \\ & & \bar{Z}_{2^{L-l}k}^L \end{bmatrix}^{-1} B_k^l. \quad (3)$$

It can be easily shown using (1)-(3) that $\bar{B}_k^l$ represents the partial scattering matrix between two adjacent level-$l$ subscatterers. This clearly suggests the butterfly compression for $\bar{B}_k^l$. The proposed factorization process proceeds as follows. At level $L$, $\bar{Z}_L$ in (1) is directly computed and inverted. At each level $l = L-1,...,0$, $\bar{Z}_l$ in (1) is computed and inverted using the following two steps. (i) Factorization: Compute a new $V = L-l$-level butterfly representation for each $\bar{B}_k^l$, $k = 1,...,2^l$ in (3). To this end, the product of $[\bar{Z}_i^s]^{-1}$ and corresponding butterfly blocks in $B_k^l$ and their partial updates are computed as a new $L-s$-level butterfly for each $s = L,...,l$ and $i = 2^{s-l}(k-1)+1,...,2^{s-l}k$. (ii) Inversion: Compute a compressed inverse for $\bar{Z}_k^{l-1}$, $k = 1,...,2^{l-1}$ in (2) as

$$[\bar{Z}_k^{l-1}]^{-1} = I + \bar{B} \quad (4)$$

where $\bar{B}$ is a new $V+1$-level butterfly. Equation (4) is suggested by the Sherman–Morrison–Woodbury formula for the inverse of a low rank update of the identity. We cannot prove but have experimentally verified the butterfly compressibility of $\bar{\boldsymbol{B}}_k^l$ in step (i) or $\bar{\boldsymbol{B}}$ in step (ii) (see Section V).

The inversion procedure for $\bar{\boldsymbol{Z}}_k^{l-1}$ in (4) can be further decomposed into five steps: (i) Split $\bar{\boldsymbol{Z}}_k^{l-1} - \boldsymbol{I}$ into four butterfly-compressed or zero submatrices, $\boldsymbol{B}_{ij}$, $i,j = 1,2$, size-wise matching the number of basis functions in the corresponding level-$l$ subscatterers. (ii) Compute a compressed inverse for $\boldsymbol{B}_{22}$ as $[\boldsymbol{I} + \boldsymbol{B}_{22}]^{-1} = \boldsymbol{I} + \bar{\boldsymbol{B}}_{22}$, where $\bar{\boldsymbol{B}}_{22}$ is a new $V$-level butterfly or zero matrix. (iii) Compute a new $V$-level butterfly representation $\tilde{\boldsymbol{B}}_{11}$ for the Schur complement of $\bar{\boldsymbol{B}}_{22}$

$$\tilde{\boldsymbol{B}}_{11} = \boldsymbol{B}_{11} - \boldsymbol{B}_{12}(\boldsymbol{I} + \bar{\boldsymbol{B}}_{22})\boldsymbol{B}_{21}. \tag{5}$$

(iv) Compute a compressed inverse for $\tilde{\boldsymbol{B}}_{11}$ as $[\boldsymbol{I} + \tilde{\boldsymbol{B}}_{11}]^{-1} = \boldsymbol{I} + \bar{\boldsymbol{B}}_{11}$, where $\bar{\boldsymbol{B}}_{11}$ is a new $V$-level butterfly. (v) Form the desired $V+1$-level butterfly $\bar{\boldsymbol{B}}$ for $[\bar{\boldsymbol{Z}}_k^{l-1}]^{-1} - \boldsymbol{I}$ in (4) from

$$\bar{\boldsymbol{B}} = \begin{bmatrix} \boldsymbol{I} & \\ -(\boldsymbol{I} + \bar{\boldsymbol{B}}_{22})\boldsymbol{B}_{21} & \boldsymbol{I} \end{bmatrix} \cdot \begin{bmatrix} \boldsymbol{I} + \bar{\boldsymbol{B}}_{11} & \\ & \boldsymbol{I} + \bar{\boldsymbol{B}}_{22} \end{bmatrix} \begin{bmatrix} \boldsymbol{I} & -\boldsymbol{B}_{12}(\boldsymbol{I} + \bar{\boldsymbol{B}}_{22}) \\ & \boldsymbol{I} \end{bmatrix} - \boldsymbol{I}. \tag{6}$$

Among the above-described five steps, steps (ii) and (iv) proceed by recursively performing steps (i)-(v).

The computational efficiency of the proposed factorization process relies on fast schemes for computing the butterfly representations of matrices on the left-hand sides of (3), (5) and (6). Note that these matrices and their transposes can be rapidly applied to arbitrary vectors as the RHSs in (3), (5) and (6) are composed of pre-computed butterfly-compressed blocks. A fast randomized butterfly scheme that relies on information gathered by multiplying the matrix with random vectors is described next.

## 4 Randomized Butterfly Scheme

Consider a $m \times n$ matrix $\boldsymbol{B}$ with $m \approx n \approx N/2^l$, its rows and columns respectively correspond to $2^v$ level-$l+v$ observation and source subscatterers, $v = 0,...,V(=L-l)$. Suppose the butterfly rank of $\boldsymbol{B}$ is capped by $r$. The proposed randomized scheme first constructs an auxiliary butterfly-factorized $m \times n$ matrix $\hat{\boldsymbol{B}}$

$$\hat{\boldsymbol{B}} = \hat{\boldsymbol{R}}^{V+1}\hat{\boldsymbol{R}}^V \cdots \hat{\boldsymbol{R}}^1\hat{\boldsymbol{R}}^0. \tag{1}$$

The diagonal blocks in $\hat{\boldsymbol{R}}^{V+1}$ have column and approximate row dimensions $r$ and $m/2^V$; similarly, the diagonal blocks in $\hat{\boldsymbol{R}}^0$ have row and approximate column dimensions $r$ and $n/2^V$; the diagonal blocks in $\boldsymbol{D}^v\hat{\boldsymbol{R}}^v$, $v = 1,...,V$ have dimensions $2r \times 2r$. All blocks in $\hat{\boldsymbol{B}}$ are filled with independent and identically distributed (i.i.d.) standard Gaussian random variables. Let $v_m = \lfloor V/2 \rfloor$ with $\lfloor \cdot \rfloor$ rounding downwards. The proposed scheme constructs $\boldsymbol{R}^{v_m} \cdots \boldsymbol{R}^0$ and $\boldsymbol{R}^V \cdots \boldsymbol{R}^{v_m+1}$ by right and left multiplying $\boldsymbol{B}$ by structured random matrices, respectively. (i) For each $v = 0,...,v_m$ and $i = 1,...,2^v$, construct a $p \times m$ structured matrix $\boldsymbol{U}$ with $p = r + O(1)$ whose columns are i.i.d. Gaussian random variables if they correspond to the $i^{th}$ level-$l+v$ observation subscatterer, and zero otherwise [see $U$ in Fig. 2(a)]. Next, compute a matrix $V'_o$ as

$$V'_o = UB(\hat{R}^0)^T \cdots (\hat{R}^{v-1})^T \qquad (2)$$

and a matrix $V'_i = V'_o(\hat{R}^v)^T$, where the superscript $T$ denotes the transpose. Note that $V'_o = UB$ when $v = 0$. It is easily shown that there are $2^{V-v}$ blocks $R$ of dimensions $r \times 2r$ [or $r \times (n/2^V)$] in $R^v$ associated with the $i^{th}$ level-$l+v$ observation subscatterer [Fig. 2(a)]. For each $R$, extract a $p \times r$ submatrix $V_i$ and a $p \times 2r$ [or $p \times (n/2^V)$] submatrix $V_o$ from $V'_i$ and $V'_o$ corresponding to the rows and columns of $R$, respectively. The block $R$ can be then computed as $R = V_i^\dagger V_o$ where † denotes the pseudoinverse. (ii) For each $v = V+1, ..., v_m+1$ and $i = 1, ..., 2^{V-v}$, construct a $n \times p$ structured matrix $U$ whose rows are i.i.d. Gaussian random variables if they correspond to the $i^{th}$ level-$L-v$ source subscatterer, and zero otherwise [Fig. 2(b)]. Compute a matrix $V'_o$

$$V'_o = (\hat{R}^{v+1})^T \cdots (\hat{R}^{V+1})^T BU. \qquad (3)$$

Furthermore, compute a matrix $V'_i = R^{v_m} \cdots R^0 U$ if $v = v_m$ and $V'_i = (\hat{R}^v)^T V'_0$ otherwise. Note that there are $2^v$ blocks $R$ of dimensions $2r \times r$ [or $(m/2^V) \times r$] in $R^v$ associated with the $i^{th}$ level-$L-v$ source subscatterer [Fig. 2(b)]. For each $R$, extract a $r \times p$ submatrix $V_i$ and a $2r \times p$ [or $(m/2^V) \times p$] submatrix $V_o$ from $V'_i$ and $V'_0$ corresponding to the columns and rows of $R$, respectively. The block $R$ can be computed as $R = V_o V_i^\dagger$. Upon completion of (i) and (ii), we have constructed a butterfly factorization $B = R^{V+1} \cdots R^0$ with blocks in $R^V, ..., R^1$ of the same dimensions. The memory of this factorization can be further reduced via applying an additional LR-compression step to all computed blocks $R$.

It can be shown that the CPU and memory costs of the above-described randomized scheme scale as $O(r^3 n^{1.5} \log n)$ and $O(r^2 n \log n)$, respectively. More importantly, the randomized scheme permits accurate construction of the butterfly with overwhelmingly high probabilities irrespective of butterfly level provided that $r$ exceeds the butterfly rank of $B$. Interestingly, we observed that the proposed direct solver can achieve good accuracy when $r$ is chosen as the maximum butterfly rank among all blocks in the RHSs of (3), (5) and (6). As a result, the proposed butterfly-based direct solver typically requires $O(N^{1.5} \log N)$ CPU and $O(N \log^2 N)$ memory resources when applied to the analysis of electrically large scatterers.

## 5 Numerical Results

This section demonstrates the applicability and efficiency of the proposed direct solver via its application to three scatterers: a corrugated semi-circle, a corrugated corner reflector, and an open cavity. For all structures, the EFIE is discretized with $0.05\lambda_0$-wide pulse basis and delta testing functions. The accuracy of the butterfly compressions is set to $10^{-4}$. All simulations are performed on a single 2.60 GHz Intel Xeon E5-2670 processor which accesses 64 GB memory.

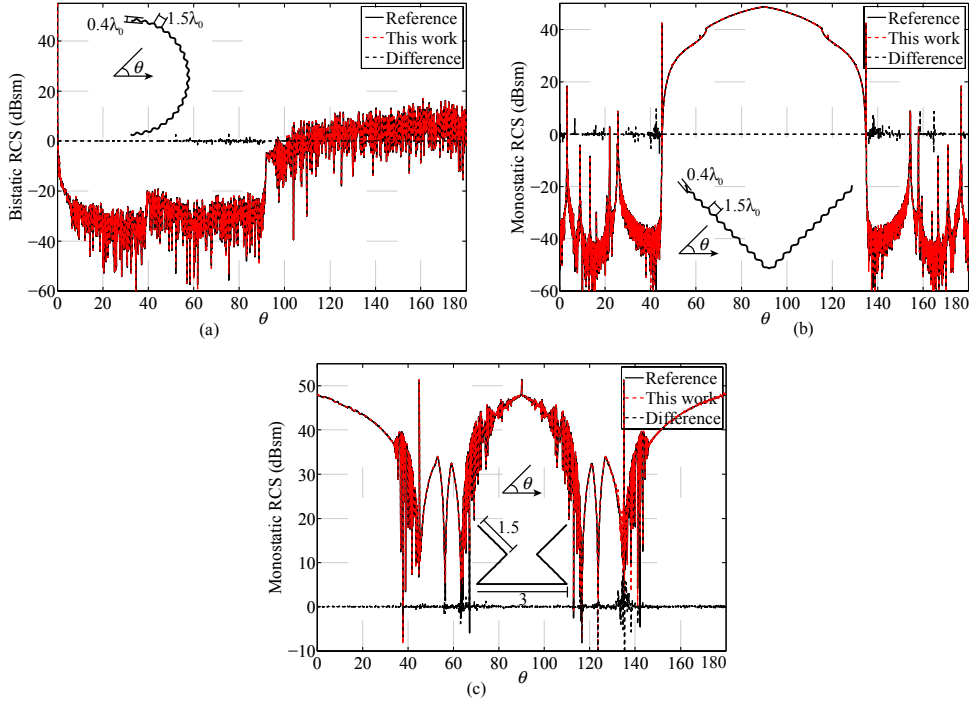

Figure 3: (a) Bistatic RCS of the corrugated semicircle. (b) Monostatic RCS of the corrugated corner reflector. (c) Monostatic RCS of the open cavity.

First, the solver is applied to a corrugated semi-circle of radius $22,084\lambda_0$. The periodicity and depth of the sinusoidal corrugations are $1.5\lambda_0$ and $0.4\lambda_0$, respectively. The structure is illuminated from the $\theta = 0^o$ direction. The corrugated semicircle is discretized with $N = 2,560,000$ basis functions and the discretized IE operator is compressed with $L = 15$ levels. The solver requires peak memory of 22.1 GB and total CPU time of 25 h. The bistatic radar cross sections (RCS) at $\theta = [0,180°]$ computed using the proposed direct solver and its LU-based counterpart [5] are in good agreement (Fig. 3(a)).

Next, the proposed solver is applied to the analysis of the monostatic RCS of the corrugated corner reflector with length $45,248\lambda_0$. Corrugation profiles, periodicity, and depth are the same as in the previous example. The structure is discretized with $N = 2,560,000$ basis functions and the discretized IE operator again compressed with $L = 15$ levels. The solver requires peak memory of 20 GB and total CPU time of 19.6 h. The monostatic RCS for 4,000 angles is computed with the proposed direct solver and that from [5] and results again agree well (Fig. 3(b)).

Next, the proposed solver is applied to the analysis of the monostatic RCS from the open cavity involving $N = 640,000$ basis functions. The solver requires the peak memory of 9 GB and total CPU time of 5.2 h. In comparison, the solver from [5] requires CPU time of 20 h. The monostatic RCS for 4,000 angles is computed using the results obtained by the proposed direct solver and the one from [5] agree well (Fig. 3(c)).

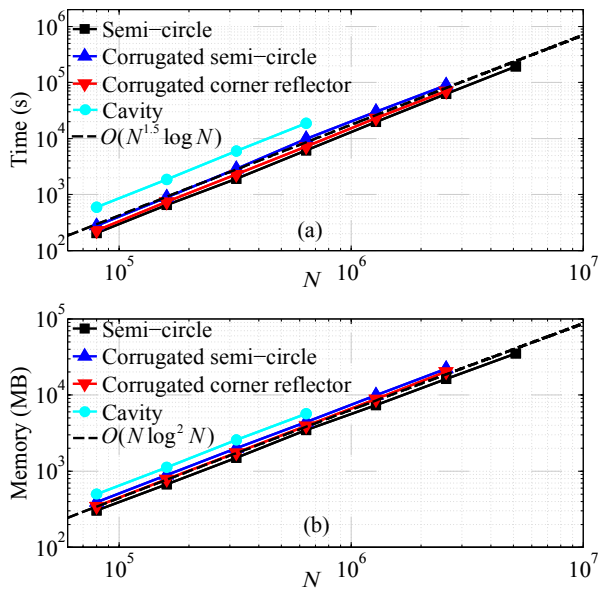

Figure 4: (a) Factorization time and (b) memory requirement of the proposed solver.

Finally, the CPU and memory resources required for inverting $Z$ when applied to the above three scatterers plus a smooth semi-circle are compared as $N$ changes from $80,000$ to $5,120,000$ (Fig. 4). The maximum butterfly ranks in $Z$ and its inverse for the smooth and corrugated semicircles, corrugated corner reflector, and cavity are 12, 20, 17 and 45, respectively. Hence the cavity requires most CPU and memory resources (per unknown). As the (observed) butterfly ranks in $Z^{-1}$ stay approximately constant irrespective the size of the scatterers, the observed CPU and memory requirements scale as $O(N^{1.5} \log N)$ and $O(N \log^2 N)$ as predicted.

## 6  Conclusion

A HSS matrix-inspired butterfly-based fast direct EFIE solver for analyzing scattering from 2D objects is presented. The solver permits butterfly compression of blocks in discretized forward and inverse EFIE operators representing near-field interactions and hinges on a fast randomized butterfly scheme. The proposed solver constitutes significant improvement over its predecessors for 2D scattering problems.

## References

[1]  S. Borm, L. Grasedyck, and W. Hackbusch, "Hierarchical matrices," *Lecture notes,* vol. 21, p. 2003, 2003.

[2] M. Bebendorf, "Hierarchical LU decomposition-based preconditioners for BEM," *Computing,* vol. 74, pp. 225-247, 2005.

[3] P.-G. Martinsson and V. Rokhlin, "A fast direct solver for boundary integral equations in two dimensions," *J. Comput. Phys.,* vol. 205, pp. 1-23, 2005.

[4] E. Corona, P.-G. Martinsson, and D. Zorin, "An O(N) direct solver for integral equations on the plane," *Appl. Comput. Harmon. Anal.,* vol. 38, pp. 284-317, 2015.

[5] H. Guo, Y. Liu, J. Hu, and E. Michielssen, "A parallel MLMDA-based direct integral equation solver," in *Proc. IEEE Int. Symp. AP-S/URSI,* 2013.

[6] E. Michielssen and A. Boag, "A multilevel matrix decomposition algorithm for analyzing scattering from large structures," *IEEE Trans. Antennas Propag.,* vol. 44, pp. 1086-1093, 1996.

[7] M. Tygert, "Fast algorithms for spherical harmonic expansions, III," *J. Comput. Phys.,* vol. 229, pp. 6181-6192, 2010.

[8] S. Ambikasaran and E. Darve, "An O(Nlog N) fast direct solver for partial hierarchically semi-separable matrices," *J. Sci. Comput.,* vol. 57, pp. 477-501, 2013.